\documentclass[11pt,reqno]{amsart}
\usepackage{amscd,amsmath,amsopn,amssymb,amsthm,multicol}
\usepackage{hyperref}
\usepackage[all]{xy}
\usepackage{amscd}
\usepackage{graphicx}

\DeclareMathOperator{\Ad}{Ad}
\DeclareMathOperator{\ad}{ad}
\DeclareMathOperator{\Aut}{Aut}
\DeclareMathOperator{\End}{End}

\DeclareMathOperator{\Iso}{Iso}

\DeclareMathOperator{\tr}{tr}
\DeclareMathOperator{\Hgt}{ht}

\DeclareMathOperator{\Ric}{Ric}
\DeclareMathOperator{\ric}{ric}

\DeclareMathOperator{\rnk}{rk}
\newcommand{\fr}{\mathfrak}
\newcommand{\al}{\alpha}
\newcommand{\be}{\beta}

\newcommand{\bb}{\mathbb}
\newcommand{\thickline}{\noalign{\hrule height 1pt}}

\newtheorem{theorem}{Theorem}
\newtheorem{lemma}{Lemma}
\newtheorem{remark}{Remark}
\newtheorem{prop}{Proposition}
\newtheorem{definition}{Definition}

\begin{document}

\title
{The Ricci flow approach to homogeneous Einstein metrics on flag manifolds}
\author{Stavros Anastassiou and Ioannis Chrysikos}
\address{Center of Research and Applications of Nonlinear Systems (CRANS), 
 University of Patras, Department of Mathematics, GR-26500 Rion, Greece}
 \email{SAnastassiou@gmail.com}
 \address{University of Patras, Department of Mathematics, GR-26500 Rion, Greece}
\email{xrysikos@master.math.upatras.gr}
  
 \medskip 
\begin{abstract}
We give the global picture of the normalized Ricci flow on  generalized  flag manifolds with two or three isotropy summands.  The normalized Ricci flow for these spaces descents to a  
parameter depending  system of two  or three ordinary differential equations, respectively.  We present here 
the qualitative study of these system's global phase portrait, by using techniques of Dynamical Systems theory.  This study  allows us to draw conclusions about the existence and  the analytical form of invariant Einstein metrics on such manifolds, and seems to offer a better insight to the classification problem of invariant Einstein metrics on compact homogeneous spaces.



\end{abstract}\maketitle

\section*{Introduction}
\markboth{Ioannis Chrysikos and Stavros Anastassiou}{The Ricci flow approach to homogeneous Einstein metrics on flag manifolds}


    The Ricci flow equation was introduced by Richard Hamilton (\cite{Ham1}) in 1982,   in order to  study  Thurston's geometrization 
conjecture on 3-manifolds. It is defined as:  
\begin{equation}\label{rf}
\frac{\partial g}{\partial t} =-2\Ric_{g}
\end{equation}
with initial condition $g(0)=g_{0}$. Here   $g=g(t)$ is a curve on the space of Riemannian metrics $\mathcal{M}$ on a smooth 
manifold $M^{n}$, and $\Ric_{g}$ is the Ricci tensor of   the Riemannian metric $g$.  Ricci flow is a way to deform a Riemannian manifold $(M, g)$
under a nonlinear evolution equation in order to process, 
or improve it.  In a more general context, Ricci flow is  a very useful tool in   proving topological 
theorems  in Riemannian 
and K\"ahlerian manifolds. The most  astonishing application  of it is
the proof of the Poincar\'e 
conjecture, due to Grisha Perelman's developments on the Hamilton's 
program.  

  Note that equation (\ref{rf}) does not preserve volume in general.
 When the underlying manifold 
$M^{n}$ is compact one often considers the {\it normalized Ricci flow} 
\begin{equation}\label{nrf}
\frac{\partial g}{\partial t}=-2\Ric_{g}+\frac{2r}{n}g,
\end{equation}
where $r=r(g(t))=\int_{M} S_{g} du_{g}/\int_{M} du_{g}$, $du_{g}$ is the volume element of $g$,   and $S_{g}$ denotes the scalar curvature function of $g$, i.e. the trace of the Ricci tensor $\Ric_{g}$.
 Under this normalized flow, the volume of the solution metric
is constant in time. Equations (\ref{rf}) and (\ref{nrf}) are 
equivalent by reparametrizing   time $t$ and scaling the metric 
in space by a function of $t$.

 Fixed points 
of the normalized Ricci flow (\ref{nrf}) are precisely the 
Einstein metrics, i.e. Riemannian  metrics  of constant 
 Ricci curvature, that is $\Ric_{g}= c\cdot g$ for some constant $c\in\mathbb{R}$. The last equation reduces to   a system of 
 second order PDEs and general existence results are difficult to    obtain. On    compact manifolds $M^{n}$ $(n\geq 3)$, Einstein metrics   of volume 1
  are in a natural way privileged metrics, since they arise as the critical 
  points of the  total scalar functional $\bold{S}(g) =\displaystyle\int_{M}S_{g}du_{g}$,
    restricted to the space $\mathcal{M}_{1}$ of all Riemannian metrics of volume 1 (cf. \cite{Be}, \cite{Wa2}).

   One of the main problems about the Ricci flow is to determine the evolution process of 
     the metrics defined on a manifold. In general, one starts with a metric $g_{0}$ on $M$ that satisfies some rather
 general curvature condition $\mathcal{R}$ and 
 proves that as the normalized Ricci flow
 runs, the metric $g_{t}$ converges to a limiting metric (as $t\rightarrow \infty$) which satisfies a more
 restrictive curvature condition $\mathcal{R}'$ (cf. \cite{Ham1},\cite{Ham2}).  
   
  Recall   that a Riemannian manifold $(M, g)$ is called $G$-homogeneous if there 
 is a closed subgroup $G$ of the group of isometries $\Iso(M, g)$  
 which acts transitively on $M$, that is, for any  $p, q\in M$  there exists $g\in G$ with $gp=q$.   
Then $M=G/K$, where   $K=\{g\in G: gp=p\}$  is the isotropy group at the point $p\in M$.   Since $G$ is a closed subgroup of $\Iso(M, g)$,  the isotropy subgroup $K$  is a compact subgroup of $\Iso(M, g)$, and thus  $M$ is compact, if and only if $G$ is compact.  When $M=G/K$ is a  homogeneous space it is 
 convenient to work with   $G$-invariant Riemannian 
 metrics $g=g(t)$, i.e. metrics for which the translations $\tau_{a} :G/K\to G/K$, $gK\mapsto agK$ act  as isometries.  For such metrics the  non-linear system  of PDEs which 
 forms the Ricci flow equation,  reduces to a 
 non-linear system  of ODEs. For this reason  one can proceed  to a 
study of the global behaviour of the Ricci flow,  using tools from 
the theory of dynamical systems (cf.  \cite{Bohm}, \cite{Glick}
 and \cite{Grama}).  
  
Let $G$ be a compact, connected and simple Lie group.  In this paper we study the global behaviour  of the normalized Ricci flow, from a 
qualitative point of view, for a $G$-invariant Riemannian metric  on   
generalized flag manifolds $M=G/K$,  whose isotropy representation decomposes into two or three pairwise inequivalent irreducible $K$-modules (see Tables 1 and 3).    Recall that a flag manifold is an adjoint orbit of $G$ and  any such space can be expressed as a   homogeneous space of the form $G/C(S)$, where $C(S)$ is the centralizer of a torus $S\subset G$.  Flag manifolds exhaust the compact and simply connected  homogeneous 
  K\"ahler manifolds  and thus they 
   have important applications in the physics of  elementary particles,
   where they give rise to a broad class  of supersymmetric sigma models.   Such a homogeneous space    admits  a finite number of $G$-invariant complex structures 
   and the first Chern class is positive. 
 Also, for any complex structure
   there is a (unique)  $G$-invariant K\"ahler-Einstein metric (cf. \cite{AP}, \cite{Arv}). 
        
The explicit form of the normalized Ricci flow is given here and its qualitative properties are presented. We determine 
the asymptotical limit for a $G$-invariant Riemannian metric of $M=G/K$ and use the results to classify all homogeneous Einstein metrics on M.
These metrics correspond to the singularities of the normalized Ricci flow   located at infinity. In particular, by using the compactificion method 
of Poincar\'{e} (\cite{Gon}), we are able to obtain the  fixed points coordinates and determine explicitly the coefficients of the Einstein metrics. Thus we prove 
the following theorems:
    
    \smallskip 
    { \bf{Theorem 1.}}
    {\it Let $M=G/K$ be a generalized flag manifold with two isotropy summands.  The normalized Ricci flow  on the space of $G$-invariant Riemmanian metrics on $M$, possesses no  finite singularities and exactly 
two  at infinity. One of them is a reppeling node and the second one is an attractive node. These fixed points determine explicity the two (up to scale) invariant Einstein metrics of $M$.}
    
 \smallskip 
    { \bf{Theorem 2.}}
    {\it Let $M=G/K$ be a generalized flag manifold of a compact simple Lie group 
$G$ with three isotropy summands, which is defined by painting black only one simple root in the Dynkin diagram of $G$.  The normalized Ricci flow  on the space of invariant Riemmanian metrics on $M$, possesses  no  finite singularities and exactly 
three   at infinity. One of them is a reppeling node while the other two are saddle points. These 
fixed points determine explicitly  the three (up to scale) $G$-invariant Einstein metrics of  $M$.}

\smallskip
The paper is structured as follows: In Section 1 we recall some basic facts of the geometry of a generalized flag manifold $M=G/K$ and we describe the form of the normalized Ricci flow for a $G$-invariant metric. In Section 2 we study flag manifolds with two or three isotropy summands and we give explicitly the system which determines  this flow. In the next section we describe the Poincar\'{e} compactification method in the two dimensional case, and  give the necessary formulas for the calculations that follow.   The final section contains the global structure of the normalized Ricci flow and the proofs of Theorems 1 and 2.

\markboth{Stavros Anastassiou and  Ioannis Chrysikos}{The Ricci flow approach to homogeneous Einstein metrics on flag manifolds}
\section{The normalized Ricci flow on generalized flag manifolds}
\markboth{Stavros Anastassiou and  Ioannis Chrysikos}{The Ricci flow approach to homogeneous Einstein metrics on flag manifolds}

\subsection{Generalized flag manifolds}
A generalized flag manifold  is an adjoint orbit  of a compact  
 semisimple Lie group  $G$.
 Recall that  $G$ acts on its Lie algebra  $\fr{g}=T_{e}G$ through 
 the adjoint action, i.e. the action induced by the adjoint representation 
 $\Ad : G\to \Aut(\fr{g})$   of $G$.  Given an element $w\in\fr{g}$ (i.e., a left-invariant vector field), 
 the adjoint orbit of $w$ is given by $M=\Ad(G)w=\{\Ad(G)w: g\in G\}\subset\fr{g}$. 
 Thus $M$ is an imbedded manifold in an Euclidean space, the Lie algebra of $G$.   Let $K=\{g\in G : \Ad(g)w=w\}\subset G$ be the isotropy subgroup of $w$ and let $\fr{k}=T_{e}K$ be the corresponding Lie algebra. 
 Since $G$ acts on $M$ transitively, $M$ is diffeomorphic to the (compact) 
 homogeneous space $G/K$, that is $\Ad(G)w=G/K$. 
  Then  one can prove  that the Lie algebra $\fr{k}$ is given by $\fr{k}=\{X\in\fr{g} : [X, w]=0\}=\ker\ad(w)$, where  $\ad :\fr{g}\to\End(\fr{g})$ is the adjoint representation of $\fr{g}$. Moreover, the set $S_{w}=\overline{\{\exp(tw) : t\in\mathbb{R}\}}$  
  is a torus in $G$ and the isotropy subgroup $K$ is 
 identified with the centralizer of  $S_{w}$ in $G$, i.e. $K=C(S_{w})$ (cf. \cite{Be}). 
 From this fact it follows   that $\rnk G=\rnk K$ and that $K$ is connected.
  
 \begin{definition}\label{def1}
 Let $G$ be  a compact semisimple Lie group. A generalized flag 
  manifold is the adjoint orbit of an element in the Lie algebra 
  $\fr{g}$ of   $G$. Equivalently, it is a homogeneous space  of 
  the form $G/K$, where $K=C(S)=\{g\in G : ghg^{-1}=h  \ \mbox{for all} \  h\in S\}$ 
  is the centralizer of a torus $S$ in $G$.  
 \end{definition}

 \subsection{Geometry of flag manifolds}
 In order to construct the normalized Ricci flow on a generalized flag manifold  it is convenient to  view it as a homogeneous space and not as an imbedded submanifold of an Euclidean space. 
 
  Let $M=G/K=G/C(S)$ be a generalized flag manifold.     Since $G$ is  compact and semisimple, the Killing form $B$
of $\fr{g}=T_{e}G$ is non-degenarate and negative definite.  
Thus the billinear form  $(   \cdot , \cdot )=-B( \cdot , \cdot )$ is an $\Ad(G)$-invariant 
 inner product on $\fr{g}$.  Let $\fr{g}=\fr{k}\oplus\fr{m}$ be a reductive 
 decomposition of $\fr{g}$ with respect to $( \cdot  , \cdot )$, that is
 $\Ad(K)\fr{m}\subset\fr{m}$.   Then,  as usual, we can identify the  $\Ad(K)$-invariant subspace $\fr{m}$ with the tangent  space $T_{o}G/K$, where $o$ denotes the identity coset of $G/K$  (cf. \cite{Arv}, \cite{Be}).   
     Under this identification, the isotropy representation  $\chi : K \to \Aut(\fr{m})$ of $G/K$ (or simply $K$)   is equivalent to the adjoint representation   
  $\Ad\big|_{K}$ restricted on $\fr{m}$, i.e. $\chi(k)=\Ad(k)\big|_{\fr{m}}$ for all $k\in K$.  Thus the
    set of  all 
  $G$-invariant symmetric  covariant  2-tensors on $G/K$  
is identified with the set of all $\Ad(K)$-invariant symmetric 
bilinear forms on $\fr{m}$.  
In particular, the set of  $G$-invariant Riemannian metrics on $G/K$ 
is identified with the set of $\Ad(K)$-invariant inner products 
$\langle \ , \ \rangle$ on $\fr{m}$ (cf. \cite{Arv}, \cite{Be}).

 Let now $g$ be a $G$-invariant metric on $M=G/K$. Since $G$ is compact and semisimple, the Ricci tensor $\Ric_{g}$ of $g$ is given by (\cite[7.38]{Be}): 
\begin{equation}\label{Ricci}
\Ric_{g}(X, X)=-\frac{1}{2}\sum_{i}|[X, X_{i}]_{\fr{m}}|^{2}+\frac{1}{2}(X, X)+\frac{1}{4}\sum_{i, j}\langle [X_{i}, X_{j}]_{\fr{m}}, X\rangle^{2},
\end{equation}
where $\{X_{i}\}$ is an orthonormal basis of $\fr{m}=T_{o}G/K$ with respect to 
the $\Ad(K)$-invariant inner product $\langle \ , \ \rangle$, induced by the $G$-invariant metric $g$.  Recall that we  can define an $\Ad(K)$-equivariant, $g$-selfadjoint endomorphism $\ric_{g}$, through the relation 
$
\Ric_{g}\big(\cdot , \cdot\big)=g\big(\ric_{g}( \ \cdot \ ) , \cdot\big).
$
 This operator is known as the Ricci operator corresponding to $g$.  Finaly, the scalar curvature $S_{g}=\tr\Ric_{g}$ of  $g$ is given by (\cite[7.39]{Be}):
\begin{equation}\label{Scale}
S_{g}=\frac{1}{2}\sum_{i}(X_{i}, X_{i})-\frac{1}{4}\sum_{i, j}|[X_{i},  X_{j}]_{\fr{m}}|^{2}.
\end{equation}
Note that  the scalar curvature $S_{g}$ is a constant function on $M=G/K$, and so   
\[
r=\frac{\int_{M} S_{g} du_{g}}{\int_{M} du_{g}}=\frac{S_{g}\int_{M}du_{g}}{\int_{m}du_{g}}=S_{g}.
\]
Thus for a $G$-invariant metric $g$ on $M=G/K$ the normalized Ricci flow (\ref{nrf})   is equivalent to  
$
\displaystyle\frac{\partial g}{\partial t}=-2\ric_{g}+\frac{2S_{g}}{n}\cdot g,
$
  where $n=\dim G/K$ (cf. \cite{Bohm}).

We   now assume that $\fr{m}=\fr{m}_1\oplus\cdots\oplus\fr{m}_{s}$
 is a $( \ , \ )$-orthogonal decomposition of $\fr{m}$
into   irreducible  mutually inequivalent $\Ad(K)$-submodules 
$\fr{m}_{j}$,  that is  $\fr{m}_{i}\ncong\fr{m}_{i}$ (as $\Ad(K)$-representations), 
for any  $1\leq i\neq j\leq s$.  Such a decomposition always exists for any flag manifold $M=G/K$ and it is expressed in terms of  {\it $\fr{t}$-roots} (cf. \cite{AP}, \cite{Arv}, \cite{Chry3}).  
It  enables us to parametrize  the set  $\mathcal{M}^{G}$ of  $G$-invariant Riemannian metrics  $g$ on 
  $G/K$ (or equivallently the set of $\Ad(K)$-invariant inner 
  products $\langle \ , \  \rangle$ on $\fr{m}$), as follows:
\begin{equation}\label{Inva}
g= \langle \ , \  \rangle= x_1\cdot ( \ , \ )|_{\fr{m}_1}+\cdots+x_s\cdot ( \ , \ )|_{\fr{m}_s},
\end{equation} 
where   $(x_1, \ldots, x_s)\in\mathbb{R}^{s}_{+}$. These metrics are diagonal with respect to the  
   decomposition $\fr{m}=\oplus_{i=1}^{s}\fr{m}_{i}$
and since  $\fr{m}_i\ncong\fr{m}_{j}$,   Schur's Lemma  implies that any $G$-invariant  metric 
on $G/K$ is of this form (cf. \cite{Wa2}). Similarly, the Ricci tensor $\Ric_{g}$, as a
$G$-invariant symmetric  covariant  2-tensor on $G/K$,  
is identified with an $\Ad(K)$-invariant symmetric billinear 
form on $\fr{m}$ and  it is given by $\Ric_{g} =\sum_{i=1}^{s} r_i\cdot ( \ , \ )|_{\fr{m}_i}$,
where $r_{1}, \ldots, r_{s}$  are the components of the Ricci 
tensor on each $\fr{m}_{i}$. Thus the  Ricci operator $\ric_{g}$ has the form 
$
\ric_{g}=\sum_{i=1}^{s}(x_i\cdot r_i)\cdot ( \ , \ )|_{\fr{m}_i},
$
where $x_1,\ldots, x_s$ are the components of the metric tensor $g$.

   Let  now $\{e_{\al}\}$ be a $( \ , \ )$-orthonormal  basis adapted to the isotropy
   decomposition $\fr{m}=\oplus_{i=1}^{s}\fr{m}_{i}$, that is  $e_{\al}\in \fr{m}_{i}$ for some $i$, 
   and $\al<\be$ if $i<j$ (with $e_{\al}\in \fr{m}_{i}$ and $e_{\be}\in\fr{m}_{j}$).  
   Let $A_{\al\be}^{\gamma}= ([e_{\al}, e_{\be}], e_{\gamma})$, so 
   that $[e_{\al}, e_{\be}]_{\fr{m}}=\sum_{\gamma}A_{\al\be}^{\gamma}e_{\gamma}$, 
   and set   
   \begin{equation}\label{tri}
   \displaystyle\genfrac{[}{]}{0pt}{}{k}{ij}=\sum(A_{\al\be}^{\gamma})^{2}=\sum ([e_{\al}, e_{\be}], e_{\gamma})^{2},
   \end{equation}
    where the sum is taken over all 
   indices $\al, \be, \gamma$ with $e_{\al}\in \fr{m}_{i}, e_{\be}\in\fr{m}_{j}$, and $e_{\gamma}\in\fr{m}_{k}$.
    These triples are called the {\it structure constants of $G/K$} with respect to the decomposition $\fr{m}=\oplus_{i=1}^{s}\fr{m}_{i}$.  Note that  $\displaystyle\genfrac{[}{]}{0pt}{}{k}{ij}$ is independent of the $( \ , \ )$-orthonormal bases choosen for 
   $\fr{m}_{i}, \fr{m}_{j}$ and $\fr{m}_{k}$, but it depends on the choise 
   of the decomposition of $\fr{m}$ (\cite{Wa2}).   
   Also the structure constants are non-negative, that is $\displaystyle { k \brack ij}\geq 0$ and $\displaystyle { k \brack ij}=0$, if and only if $([\fr{m}_{i}, \fr{m}_{j}], \fr{m}_{k})=0$ and they are symmetric in all 
   three entries, that is $\displaystyle\genfrac{[}{]}{0pt}{}{k}{ij}=\genfrac{[}{]}{0pt}{}{k}{ji}=\genfrac{[}{]}{0pt}{}{j}{ki}$.
By reconstructing an $\langle \ , \ \rangle$-orthonormal basis for each isotropy summand $\fr{m}_{k}$ and  applying relations (\ref{Ricci}) and (\ref{Scale}), we obtain the following useful expressions.
 \begin{prop}\label{Ricc}{\textnormal{(\cite{Wa2}, \cite{SP})}} 
Let $M=G/K$ be a generalized flag manifold of a compact semisimple Lie group $G$ and let
 $\fr{m}=\oplus_{i=1}^{s}\fr{m}_{i}$ be an isotropy decomposition of    $\fr{m}$. 
  Let $g$ be 
 a $G$-invariant metric on $M$ defined by (\ref{Inva}), and set $d_{k}=\dim \fr{m}_{k}$ for all $k=1, \ldots s$. Then:
 
(1)  The components $r_{1}, \ldots, r_{s}$ of the Ricci tensor $\Ric_{g}$  are given by
   \begin{equation*} 
   r_{k}=\frac{1}{2x_{k}}+\frac{1}{4d_{k}}\sum_{i, j}\frac{x_{k}}{x_{i}x_{j}} \genfrac{[}{]}{0pt}{}{k}{ij} -\frac{1}{2d_{k}}\sum_{i, j}\frac{x_{j}}{x_{k}x_{i}} \genfrac{[}{]}{0pt}{}{j}{ki}, \qquad (k=1, \ldots, s).
 \end{equation*}
 
 (2) The scalar curvature $S_{g}={\rm tr}\Ric_{g}=\sum_{i=1}^{s}d_{i}r_{i}$ is given by
 $
  S=\displaystyle\frac{1}{2}\sum_{i=1}^{s}\frac{d_{i}}{x_{i}}-\frac{1}{4}\sum_{i, j, k} \genfrac{[}{]}{0pt}{}{k}{ij}\frac{x_{k}}{x_{i}x_{j}}.
 $
  \end{prop}
  It follows that  for a $G$-invariant metric $g=\langle \ , \  \rangle$ given by (\ref{Inva}), the normalized Ricci flow  reduces to the following system (we ommit the minus sign of the Ricci operator since our $G$-invariant metric $g=\langle \ , \  \rangle$ is expressed with respect to the inner product $( \ , \ )=-B( \ , \ )$):
 
  \begin{equation}\label{nrff}
  \left. \begin{tabular}{lcl}
  $\dot{x}_1$ & $=$ & $2x_1 \cdot r_1 + \displaystyle\frac{2 S_{g}}{n}\cdot x_1$, \\\\
$\dot{x}_2$ & $=$ &  $2x_2\cdot r_2  + \displaystyle\frac{2 S_{g}}{n}\cdot x_2,$  \\
 & $\vdots$ &   \\
$\dot{x}_s$ & $=$ & $2x_s\cdot r_s  + \displaystyle\frac{2 S_{g}}{n}\cdot x_s.$
\end{tabular}\right\}
\end{equation}

 \markboth{Stavros Anastassiou and  Ioannis Chrysikos}{The Ricci flow approach to homogeneous Einstein metrics on flag manifolds}
\section{Flag manifolds with two and three isotropy summands}
\markboth{Stavros Anastassiou and  Ioannis Chrysikos}{The Ricci flow approach to homogeneous Einstein metrics on flag manifolds}
 
For the construction of flag manifolds with two or three isotropy summands,  it is useful to describe   the structure of such a space $M=G/K=G/C(S)$    in terms of Lie theory.  By using this  description one can classify flag manifolds thgrough the {\it painted Dynkin diagrams}.    
    
 \subsection{A description of flag manifolds in terms of painted Dynkin diagrams}
     For simplicity we assume that $G$ is simple with Lie algebra $\fr{g}$. 
  Let $ \fr{g}^{\mathbb{C}}=\fr{h}^{\mathbb{C}}\oplus\sum_{\al\in R}\fr{g}_{\al}^{\mathbb{C}}$
  be the root space decomposition of the complexification 
  $\fr{g}^{\mathbb{C}}$ of $\fr{g}$, with respect to a 
  Cartan subalgebra $\fr{h}^{\mathbb{C}}$ of $\fr{g}^{\mathbb{C}}$, 
  where $R\subset(\fr{h}^{\mathbb{C}})^{*}$ is the root system of 
  $\fr{g}^{\mathbb{C}}$ and 
  $\fr{g}_{\al}^{\mathbb{C}}=\{X\in\fr{g}^{\bb{C}} : \ad(H)X=\al(H)X, \ \mbox{for all} \  H\in\fr{h}^{\bb{C}}\}$ 
  are the 1-dimensional root spaces. 
   As usual, we identify $(\fr{h}^{\mathbb{C}})^{*}$ 
   with $\fr{h}^{\mathbb{C}}$ via the Killing form $B$ 
   of $\fr{g}^{\mathbb{C}}$.   
   Let  $\Pi=\{\al_{1}, \ldots, \al_{\ell}\}$ \ 
   $(\dim\fr{h}^{\mathbb{C}}=\ell)$ be a fundamental system of 
   $R$ and choose a subset  $\Pi_{K}$   of $\Pi$. We denote by 
   $R_{K}=\{\be\in R : \be=\sum_{\al_{i}\in\Pi_{K}}k_{i}\al_{i}\}$   
   the closed subsystem spanned by $\Pi_{K}$. Then the Lie subalgebra  $  \fr{k}^{\mathbb{C}}=\fr{h}^{\mathbb{C}}\oplus\sum_{\be\in R_{K}}\fr{g}_{\be}^{\mathbb{C}}$
    is a reductive subalgebra of $\fr{g}^{\bb{C}}$, i.e. 
    it admits a decomposition of the form  $\fr{k}^{\bb{C}}=Z(\fr{k}^{\bb{C}})\oplus\fr{k}_{ss}^{\bb{C}}$,
     where $Z(\fr{k}^{\bb{C}})$ is its center and 
     $\fr{k}_{ss}^{\bb{C}}=[\fr{k}^{\bb{C}}, \fr{k}^{\bb{C}}]$ 
     the semisimple part of $\fr{k}^{\bb{C}}$. In particular, 
      $R_{K}$ is the root system of $\fr{k}^{\bb{C}}_{ss}$, and 
      thus $\Pi_{K}$ can be considered as the associated fundamental  
      system.  Let    $K$ be the connected Lie subgroup of $G$
       generated by $\fr{k}=\fr{k}^{\mathbb{C}}\cap\fr{g}$. 
        Then the homogeneous manifold $M=G/K$ is a flag manifold, 
        and any flag manifold is defined in this way, i.e., 
        by the choise of a triple $(\fr{g}^{\bb{C}}, \Pi, \Pi_{K})$ (cf. \cite{AA}, \cite{AP}).

      Set  $\Pi_{M}=\Pi\backslash \Pi_{K}$  and  $R_{M}=R\backslash R_{K}$, such that   $\Pi=\Pi_{K}\sqcup \Pi_{M}$, and $R=R_{K}\sqcup R_{M}$, respectively.  Roots in  $R_{M}$ are called  {\it complementary roots}, and they possess an important role in the geometry of $M=G/K$.  For example, if $\fr{g}^{\bb{C}}=\fr{k}^{\bb{C}}\oplus\fr{m}^{\bb{C}}$ is a reductive decomposition of $\fr{g}^{\bb{C}}$, then   $\fr{m}^{\bb{C}}=(T_{o}G/K)^{\bb{C}}=\sum_{\al\in R_{M}}\fr{g}_{\al}^{\bb{C}}$.  
 
 \begin{definition}\label{pdd}
 Let $\Gamma=\Gamma(\Pi)$ be the Dynkin diagram of the fundamental 
 system $\Pi$.  By painting in  black the nodes 
 of $\Gamma$  corresponding    to  $\Pi_{M}$, we obtain the painted 
 Dynkin diagram of the flag manifold $G/K$. In this diagram
    the subsystem $\Pi_{K}$ is determined as the subdiagram of white roots. 
 \end{definition}
  
      Conversely, given a painted Dynkin diagram, in order to determine the associated flag manifold $M=G/K$ we are working as follows: At first, we obtain the group $G$, as the unique simply connected Lie group generated by the  unique real   form $\fr{g}$   of the   complex simple Lie algebra $\fr{g}^{\bb{C}}$ (up to inner automorphisms of $\fr{g}^{\bb{C}}$ \cite[p. ~184]{Hel}),  which is reconstructed by the  underlying Dynkin diagram. On the other hand, the connected Lie subgroup $K\subset G$ is defined by using  the additional information $\Pi=\Pi_{K}\cup\Pi_{M}$ coded into the painted Dynkin diagram.  The semisimple part of $K$, is obtained from  the (not  necessarily connected) subdiagram of white roots, and   each black root, i.e.   each  root in $\Pi_{M}$,  gives rise to one $U(1)$-summand.      Thus  the painted Dynkin diagram determines 
    the  isotropy group $K$ and the space $M=G/K$ completely.        By using certain rules to determine 
     whether     different painted Dynkin diagrams define isomorphic 
    flag manifolds (see \cite{AA} , \cite{AP}),  one  can obtain all flag manifolds $G/K$ of  a 
    compact  simple Lie group $G$ .

\subsection{Flag manifolds with two isotropy summands}
Let $G$ be a compact simple Lie group. We will 
construct the normalized Ricci flow equation for  
$G$-invariant metrics on flag manifolds $M=G/K$ with two 
isotropy summands, that is $\fr{m}=\fr{m}_1\oplus\fr{m}_2$. 
 These spaces have been classified in terms of painted Dynkin 
 diagrams in \cite{Chry1} (see also \cite{Chry2}, \cite{Sak}), 
 and they are given in Table 1.  In particular, they are determined by painting black 
 in the Dynkin diagramm of $G$, a simple root $\al_{p}\in\Pi$ with height 2, that is $\Pi\backslash\Pi_{K}=\Pi_{M}=\{\al_{p} : \Hgt(\al_{p})=2\}$.\footnote{The height of a simple root $\al_{p}\in\Pi$  $(p=1, \ldots, \ell)$, is the  positive integer   $c_{p}$ in 
the expression of the highest root $\widetilde\alpha=\sum_{k=1}^{\ell}c_{k}\al_{k}$ of $\fr{g}^{\bb{C}}$, in terms of simple roots.  Note that  $\Hgt : \Pi\to\bb{Z}^{+}$  
  is the  function defined by $\Hgt(\al_{p})=c_{p}$.}

 \begin{center}
{\bf{Table 1.}} {\small The generalized flag manifolds with two isotropy summands. }
  \end{center}  
   \begin{center}
{\footnotesize{
$
   \begin{tabular}{ll}
   \hline
   $G \ \mbox{simple}$ & $\mbox{Flag manifold} \ \ G/K \  \ \mbox{with} \ \ \ \fr{m}=\fr{m}_1\oplus\fr{m}_2$ \\
   \thickline
$B_{\ell}$ &   $SO(2\ell +1)/U(p)\times SO(2(\ell-p)+1) \ \ (2\leq p\leq\ell)$ \\
$C_{\ell}$ &  $Sp(\ell)/U(p)\times Sp(\ell-p) \ \ (1\leq p\leq \ell-1)$\\
$D_{\ell}$ & $SO(2\ell)/U(p)\times SO(2(\ell-p)) \ \ (2\leq p\leq \ell-2)$   \\
$G_2$ &   $G_{2}/U(2)$ \   ($U(2)$ is represented by the short root of $G_{2}$) \\
$F_4$ &  $ F_{4}/SO(7)\times U(1)$ \\
& $F_{4}/Sp(3)\times U(1)$ \\
 $E_6$ &  $E_{6}/SU(6)\times U(1)$\\
 &    $E_{6}/SU(2)\times SU(5)\times U(1)$\\
$E_7$ &  $E_{7}/SU(7)\times U(1)$\\
&   $E_{7}/ SU(2)\times SO(10)\times U(1)$\\
 &  $E_{7}/ SO(12)\times U(1)$\\
 $E_8$ & $E_{8}/E_{7}\times U(1)$\\
 & $E_{8}/SO(14)\times U(1)$\\
  \hline
   \end{tabular}
  $}}
   \end{center}

\medskip
 We mention here, for future use, that the isotropy summands $\fr{m}_1$ and $\fr{m}_{2}$ satisfy the following useful inclusions (see \cite{Chry2}, or \cite{Do}):
\begin{equation}\label{incl2}
[\fr{m}_1, \fr{m}_1]\subset\fr{k}\oplus\fr{m}_2, \quad [\fr{m}_1, \fr{m}_2]\subset\fr{m}_1, \quad [\fr{m}_2, \fr{m}_2]\subset\fr{k}.
\end{equation} 
Now, according to  (\ref{Inva}), a $G$-invariant metric $g$ of $M$ is determined by two positive variables, i.e.
\begin{equation}\label{metric1}
g=\langle \ , \ \rangle=x_1 \cdot ( \ , \ )|_{\fr{m}_1}+x_2\cdot ( \ , \ )|_{\fr{m}_2}, \quad  x_1>0, \ x_2>0. 
\end{equation}
Thus the space $\mathcal{M}^{G}$ of $G$-invariant metrics on $G/K$ is 2-dimensional (since $\fr{m}_{1}\ncong\fr{m}_{2}$).  From inclusions (\ref{incl2}) and by applying relation (\ref{tri}), we easily conlude that the only non-zero structure constants $\displaystyle {k \brack {ji}}$ of $G/K$ are the triples
$
 \displaystyle{2 \brack 11}={1 \brack 12}={1 \brack 21}\neq 0.
$
 Now the Ricci tensor $\Ric_{g}$ of $(M, g)$ is given by 
 $
 \Ric_{g}=r_1\cdot ( \ , \ )|_{\fr{m}_1}+r_2\cdot ( \ , \ )|_{\fr{m}_2},
 $
 where  $r_1, r_2\in\bb{R}$.   
  Thus by applying Proposition  \ref{Ricc}  we easily obtain 
   \begin{prop}\label{info2}
 Let $M=G/K$ be a generalized flag manifold with two isotropy summands and let $g$ be a  $G$-invariant Riemannian metric on $M$ given by (\ref{metric1}). 
 
 (1) The components $r_1, r_2$ of the Ricci tensor $\Ric_{g}$  are given as follows:
\[
\left\{ 
\begin{tabular}{l}
$r_1=\displaystyle\frac{1}{2x_1} - \displaystyle{{2 \brack 11}}\displaystyle\frac{x_2}{{2d_1}{x_{1}^2}},$ \\\\ 
$r_2=\displaystyle\frac{1}{2x_2} +  \displaystyle{{2 \brack 11}}\Big(\displaystyle\frac{x_2}{{4d_2 x_{1}^2}}  -\displaystyle\frac{1}{2d_2x_2}\Big)$.
\end{tabular}\right.
\]

(2) Ths scalar curvature $S_{g}$ is given by 
$
S_{g}=\sum_{i=1}^{2}d_{i}\cdot r_{i}=\displaystyle\frac{1}{2}\Big(\frac{d_{1}}{x_{1}}+\frac{d_{2}}{x_{2}}\Big)-\frac{1}{4}\displaystyle{{2 \brack 11}}\Big(\frac{x_{2}}{x_{1}^{2}}+2\frac{1}{x_{2}}\Big).
$
 \end{prop}

 In order to find the non-zero structure constant $\displaystyle{{2 \brack 11}}$ of $M=G/K$ (with respect the decomposition   $\fr{m}=\fr{m}_1\oplus\fr{m}_{2}$), we use the unique $G$-invariant K\"ahler-Einstein metric $g_{J}$ which admits $M=G/K$, compatible with the unique $G$-invariant complex structure $J$ on $M$ (cf. \cite[13.8]{Borel}).  This invariant metric is given by  (cf. \cite{Chry2})
 \[
 g_{J}=1\cdot ( \ , \ )|_{\fr{m}_1}+2\cdot ( \ , \ )|_{\fr{m}_2}.
 \]
  Substituting  the values   $x_1=1$, and $x_2=2$ in  $r_1-r_2=0$, we easily obtain that  (cf.  \cite{Chry2})
\[
\displaystyle{{2 \brack 11}=\frac{d_1d_2}{d_1+4d_2}}.
\]  
The dimensions $d_{i}=\dim\fr{m}_{i}$ were calculated in \cite{Chry2} and   are given in the following table:

 \begin{center}
{\bf{Table 2.}} {\small The dimensions $d_{i}=\dim\fr{m}_{i}$ for any $M=G/K$ with $\fr{m}=\fr{m}_1\oplus\fr{m}_{2}$}
  \end{center}  

 \begin{center} 
    {\footnotesize
        \begin{tabular}{l l l}
 \hline\hline  
     $M=G/K$       &  $d_1$ & $d_2$    \\
  \hline\hline
$SO(2\ell +1)/U(p)\times SO(2(\ell-p)+1) \ \ (2\leq p\leq\ell)$                 &  $2p(2(\ell-p)+1)$  & $p(p-1)$ \\
$Sp(\ell)/U(p)\times Sp(\ell-p) \ \ (1\leq p\leq \ell-1)$                       &  $4p(\ell-p)$       & $p(p+1)$     \\
$SO(2\ell)/U(p)\times SO(2(\ell-p)) \ \ (2\leq p\leq \ell-2)$                   &  $4p(\ell-p)$       & $p(p-1)$   \\ 
$G_{2}/U(2) $                                                                      &  $8$   & $2$   \\
$F_{4}/SO(7)\times U(1)$                                                           &  $16$  & $14$  \\
$F_{4}/Sp(3)\times U(1)$                                                           &  $28$  & $2$   \\ 
$E_{6}/SU(6)\times U(1)$                                                           &  $40$  & $2$     \\
$E_{6}/SU(2)\times SU(5)\times U(1)$                                               &  $40$  & $10$    \\
$E_{7}/SU(7)\times U(1)$                                                           &  $70$  & $14$     \\
$E_{7}/ SU(2)\times SO(10)\times U(1)$                                             &  $64$  & $20$  \\
$E_{7}/ SO(12)\times U(1)$                                                         &  $64$  & $2$    \\
$E_{8}/E_{7}\times U(1)$                                                           &  $112$ & $2$      \\
$E_{8}/SO(14)\times U(1)$                                                          &  $128$ & $28$     \\
\hline
    \end{tabular}}
 \end{center}
Now,  the Ricci components $r_1, r_2$ and the scalar curvature $S_{g}$  are completely determined from Proposition \ref{info2}. As a consequence of relation (\ref{nrff})   we obtain that the normalized Ricci flow for a homogeneous initial metric (\ref{metric1}) on $M=G/K$ is given by  
 \begin{equation}\label{nrf2}
\left.\begin{tabular}{ll}
$\dot{x}_1 =$ & $2x_1 \cdot r_1 + \displaystyle\frac{2x_1}{d_1 + d_2}\cdot S_{g}$, \\
$\dot{x}_2 =$ & $2x_2\cdot r_2  + \displaystyle\frac{2x_2}{d_1 + d_2}\cdot S_{g}$.
\end{tabular}\right\}
\end{equation}

 \begin{remark}\label{Einstein1}
 \textnormal{In \cite{Ker} Dickinson and Kerr   applied the variational method to show  that the number of $G$-invariant Einstein metrics on $M=G/K$ is two.  The explicit form of these  metrics (which are not isometric) was given in \cite{Chry2} by solving the equation $r_1-r_2=0$, substituting first the above value of the triple $\displaystyle {2 \brack 11}$.  The first one is the above defined  K\"ahler-Einstein metric, and the second one   is non-K\"ahler and it is given  by 
 \[
 g=1\cdot( \ , \ )|_{\fr{m}_1}+\displaystyle\frac{4d_2}{d_1+2d_2}\cdot ( \ , \ )|_{\fr{m}_2}.
 \]
}
\end{remark}

\subsection{Flag manifolds with three isotropy summands}
Flag manifolds $M=G/K$ of a compact simple Lie group $G$, whose isotropy representation decomposes into three pairwise non isomorphic irreducible $\Ad(K)$-modules, i.e. $\fr{m}=\fr{m}_1\oplus\fr{m}_2\oplus\fr{m}_3$,
 are obtained by painting the Dynkin diagram $\Gamma(\Pi)$ of $G$ with two differents ways. We can either paint a simple root with height 3, or two simple roots both of height 1.  Thus the pairs $(\Pi, \Pi_{K})$ which determine such spaces are devided into two different types  as follows (cf. \cite{Kim}):
  \medskip
\begin{center}
 \begin{tabular}{r|l}
     Type  & $(\Pi, \Pi_{K})$\\
 \hline  
 I & $ \Pi\setminus\Pi_{K}=\{\al_{p} :  \Hgt(\al_{p})=3\}$ \\
 \hline 
  II & $\Pi\setminus\Pi_{K}=\{\al_{p}, \al_{q} \ (1\leq p\neq q\leq\ell) : \Hgt(\al_{p})=\Hgt(\al_{q})=1\}$\\ 
 \hline
  \end{tabular}
 \end{center}
  
  \medskip
  In the following we shall call the  flag manifolds   determined by a pair $(\Pi, \Pi_{K})$ of Type I (resp. of Type II), {\it flag manifolds of Type I} (resp. {\it flag manifolds of Type II}). We present these homogeneous spaces in Table 3.  

 \begin{center}
{\bf{Table 3.}} {\small The generalized flag manifolds with three isotropy summands. }
  \end{center}  
   
{\small \begin{center}
   \begin{tabular}{ll}
   \hline
  $G \ \mbox{simple}$ & $\mbox{Flag manifold} \ \ G/K  \ \ \mbox{of Type  I}$ \\
   \thickline
$E_8$ & $E_8/E_6\times SU(2)\times U(1)$  \\ 
&  $E_8/SU(8)\times U(1)$  \\
$E_7$ &  $E_7/SU(5)\times SU(3)\times U(1)$   \\
& $E_7/SU(6)\times SU(2)\times U(1)$    \\
$E_6$ & $E_6/SU(3)\times SU(3)\times SU(2)\times U(1)$ \\
  $F_4$ &  $F_4/SU(3)\times SU(2)\times U(1)$ \ ($SU(2)$ is represented by the long root of $F_4$)   \\
 $G_2$ & $G_2/U(2)$  \ ($U(2)$ is represented by the long root of $G_{2}$)   \\
  \hline 
   $G \ \mbox{simple}$ & $\mbox{Flag manifold} \ \ G/K   \ \ \mbox{of Type  II}$ \\
    \thickline 
  $A_{\ell}$ &  $SU(\ell+m+n)/S(U(\ell)\times U(m)\times U(n))$ \ ($\ell, m, n\in\bb{Z}^{+}$) \\
  $D_{\ell}$ &  $SO(2\ell)/U(1)\times U(\ell-1)$ \ ($\ell\geq 4$) \\
  $E_6$ & $E_6/SO(8)\times U(1)\times U(1)$ \\
  \hline
   \end{tabular}
        \end{center}}

  \subsection{Flag manifolds of Type I}
 Let us  now  construct  the normalized Ricci flow for a flag manifold $M=G/K$ of Type I with $\fr{m}=\fr{m}_1\oplus\fr{m}_2\oplus\fr{m}_3$.  The case of the flag manifolds of Type II  will be treated in a forthcoming paper (see also \cite{Grama} for the full flag manifold $SU(3)/T$). 
 
  The dimensions $d_{i}=\dim\fr{m}_{i}$ $(i=1, 2, 3)$ have been calculated  in \cite[p.~ 311]{Kim}, and are presented in Table 4.
  
  \begin{center}
{\bf{Table 4.}} {\small The dimensions $d_{i}=\dim\fr{m}_{i}$ for any $M=G/K$ of Type I.}
  \end{center}  
   
{\small \begin{center}
   \begin{tabular}{llll}
   \hline
   $\mbox{Flag manifold} \ \ G/K  \ \ \mbox{of Type  I}$ & $d_1$ & $d_2$ & $d_3$ \\
   \thickline
 $E_8/E_6\times SU(2)\times U(1)$  & $108$ & $54$ & $4$ \\ 
  $E_8/SU(8)\times U(1)$  & $112$ & $56$ & $16$ \\
  $E_7/SU(5)\times SU(3)\times U(1)$ & $60$ & $30$ & $8$  \\
 $E_7/SU(6)\times SU(2)\times U(1)$  & $60$ & $30$ & $4$  \\
  $E_6/SU(3)\times SU(3)\times SU(2)\times U(1)$ & $36$ & $18$ & $4$ \\
   $F_4/SU(3)\times SU(2)\times U(1)$  & $24$ & $12$ & $4$ \\
   $G_2/U(2)$   & $4$ & $2$ & $4$    \\
  \hline 
       \end{tabular}
        \end{center}}

\medskip  
 Note that the isotropy summands  $\fr{m}_{1},\fr{m}_{2}$ and $\fr{m}_{3}$ fullfil the following inclusions   (see  \cite{Do}):
 \begin{equation}\label{incl3}
 \left.\begin{tabular}{lll}
 $[\fr{m}_1, \fr{m}_1]\subset\fr{k}\oplus\fr{m}_2,$ & $[\fr{m}_1, \fr{m}_2]\subset\fr{m}_1\oplus\fr{m}_{3},$ & $[\fr{m}_1, \fr{m}_3]\subset\fr{m}_2,$ \\
 $[\fr{m}_2, \fr{m}_2]\subset\fr{k},$ & $[\fr{m}_2, \fr{m}_3]\subset\fr{m}_1,$ & $[\fr{m}_3, \fr{m}_3]\subset\fr{k}$. 
 \end{tabular} \right\}
 \end{equation}
  According to  (\ref{Inva}), a $G$-invariant metric $g$ of $M$ is determined by three  positive parameters, i.e.
\begin{equation}\label{metric2}
g=\langle \ , \ \rangle=x_1 \cdot ( \ , \ )|_{\fr{m}_1}+x_2\cdot ( \ , \ )|_{\fr{m}_2}+x_3 \cdot ( \ , \ )|_{\fr{m}_3}, \quad  x_1>0, \ x_2>0, \  x_3>0, 
\end{equation}
and since $\fr{m}_{1}\ncong\fr{m}_{2}\ncong\fr{m}_{3}$, the space $\mathcal{M}^{G}$ of $G$-invariant metrics on $G/K$ is 3-dimensional.  
   In view of   (\ref{incl3}) and  by applying relation (\ref{tri}), we conclude that the only non zero structure constants of $M=G/K$ are the following
   \begin{eqnarray}
    { 1 \brack 12}  &=&  { 1 \brack 21} = { 2 \brack 11}=c_{11}^{2} \nonumber \\
    { 3 \brack 12} &=& { 3 \brack 21} = { 2 \brack 13} = { 2 \brack 31} = { 1 \brack 23} = { 1 \brack 32}= c_{12}^{3}. 
    \end{eqnarray} 
 Now, the Ricci tensor $\Ric_{g}$ of $g$ is given by $\Ric_{g}=r_1\cdot ( \ , \ )|_{\fr{m}_1}+r_2\cdot ( \ , \ )|_{\fr{m}_2}+r_3\cdot ( \ , \ )|_{\fr{m}_3},
$ where   $r_1, r_2, r_3\in\bb{R}$.   
  Thus, by applying Proposition  \ref{Ricc}, we easily obtain the following.

 \begin{prop}\label{info}
 Let $M=G/K$ be a generalized flag manifold with three isotropy summands of Type I and let $g$ be a  $G$-invariant Riemannian metric on $M$ given by (\ref{metric2}). 
 
 (1) The components $r_1, r_2, r_3$ of the Ricci tensor $\Ric_{g}$  are given as follows:
\[
\left\{ 
\begin{tabular}{l}
$r_1=\displaystyle\frac{1}{2x_1} - \displaystyle\frac{c_{11}^{2}x_2}{2d_1x_{1}^{2}}+\displaystyle\frac{c_{12}^{3}}{2d_1}\Big(\frac{x_1}{x_2x_3}-\frac{x_2}{x_1x_3}-\frac{x_3}{x_1x_2}\Big),$ \\\\ 
$r_2=\displaystyle\frac{1}{2x_2} +  \displaystyle\frac{c_{11}^{2}}{4d_2}\Big(\frac{x_2}{x_{1}^2} -\frac{2}{x_2}\Big)+\displaystyle\frac{c_{12}^{3}}{2d_2}\Big(\frac{x_2}{x_1x_3}-\frac{x_1}{x_2x_3}-\frac{x_3}{x_1x_2}\Big),$ \\\\
$r_3=\displaystyle\frac{1}{2x_3} +  \displaystyle\frac{c_{12}^{3}}{2d_3}\Big(\frac{x_3}{x_1x_2}-\frac{x_1}{x_2x_3}-\frac{x_2}{x_1x_3}\Big).$
\end{tabular}\right.
\]

(2) Ths scalar curvature $S_{g}$  is given by 
\[ 
S_{g}=\sum_{i=1}^{3}d_{i}\cdot r_{i}=\frac{1}{2}\Big(\frac{d_{1}}{x_{1}}+\frac{d_{2}}{x_{2}}+\frac{d_{3}}{x_{3}}\Big)-\displaystyle\frac{c_{11}^{2}}{4}\Big(\frac{x_{2}}{x_{1}^{2}}+ \frac{2}{x_{2}}\Big) -  \displaystyle\frac{c_{12}^{3}}{2}\Big(\frac{x_1}{x_2x_3}+\frac{x_2}{x_1x_3}+\frac{x_3}{x_1x_2}\Big).
\]
 \end{prop}

  For the computation of the unknown triples $c_{11}^{2}$ and $c_{12}^{3}$ we use the  unique  $G$-invariant K\"ahler-Einstein metric which admits any flag manifold $G/K$ of Type I, and given by (cf. \cite{Kim}):
  \[
  g_{J}=\langle \ , \ \rangle= 1\cdot ( \ , \ )|_{\fr{m}_1}+2\cdot ( \ , \ )|_{\fr{m}_2}+3\cdot ( \ , \ )|_{\fr{m}_3}.
  \]
    Substituting the values   $x_1=1$,  $x_2=2$ and $x_3=3$ in the system
  \begin{equation}\label{e1}
   r_1-r_2=0, \quad r_2-r_3=0,
  \end{equation}
 where $r_1, r_2, r_3$ are given by Proposition \ref{info},  we  obtain the following values:
 \begin{equation}\label{stru3}
 c_{11}^{2}=\displaystyle\frac{d_1 d_2 + 2 d_1 d_3 - d_2 d_3}{d_1 + 4 d_2 + 9 d_3}, \qquad c_{12}^{3}=\displaystyle\frac{(d_1 + d_2) d_3}{d_1 + 4 d_2 + 9 d_3}.
 \end{equation}

 Thus, the Ricci components $r_1, r_2, r_3$ and the scalar curvature $S_{g}$  are completely determined from Proposition \ref{info}.  By applying relation  (\ref{nrff}), we obtain  the normalized Ricci flow for a homogeneous initial metric (\ref{metric2}) on $M=G/K$, defined  by the following system:
  \begin{equation}\label{nrf3}
 \left. \begin{tabular}{ll}
 $\dot{x}_{1} =$ & $2x_1\cdot r_1+\displaystyle\frac{2x_1}{d_1+d_2+d_3}\cdot S_{g}$, \\
 $\dot{x}_{2} =$ & $2x_2\cdot r_2+\displaystyle\frac{2x_2}{d_1+d_2+d_3}\cdot S_{g}$, \\
 $\dot{x}_{3} =$ & $2x_3\cdot r_3+\displaystyle\frac{2x_3}{d_1+d_2+d_3}\cdot S_{g}$. 
 \end{tabular}\right\}
 \end{equation}

 \begin{remark}\label{Einstein1}
 \textnormal{In \cite{Kim}, was proved that  the number of $G$-invariant Einstein metrics on a flag manifold $M=G/K$ of Type I is exactly three.  One is the above defined K\"ahler--Einstein metric $g_{J}$ and the other two are non K\"ahler.   The explicit form  however of these two metrics was not presented. We give here these explicit forms by solving  equation  (\ref{e1}). 	First we normalize the $G$-invariant metric $(\ref{metric2})$ by setting $x_1=1$.
 We obtain the following theorem:}
 \end{remark}
 
 \begin{theorem}\label{type13} \textnormal{(\cite{Kim})}
 Let $M=G/K$ be a generalized flag manifold of a compact simple Lie group 
$G$ of Type I. Then $M$ admits precisely    three (up to a scale) $G$--invariant Einstein 
metrics. One is K\"ahler--Einstein  given by $g_{J}=(1, 2, 3)$, and the other two are non K\"ahler metrics $g_1, g_2$, given   below: 
{\small \begin{center}
   \begin{tabular}{lll}
   \hline
   $\mbox{Flag manifold} \ \ G/K  \ \ \mbox{of Type  I}$ & $g_{1}=(1, x_2, x_3)$ & $g_{2}=(1, x_2, x_3)$ \\
   \thickline
 $E_8/E_6\times SU(2)\times U(1)$  & $(1, 0.914286, 1.54198)$ & $(1, 1.0049,  0.129681)$   \\ 
  $E_8/SU(8)\times U(1)$  & $(1, 0.717586,  1.25432)$ & $(1, 1.06853,  0.473177)$   \\
  $E_7/SU(5)\times SU(3)\times U(1)$ & $(1, 0.733552,  1.27681)$ & $(1, 1.06029, 0.443559)$   \\
 $E_7/SU(6)\times SU(2)\times U(1)$  & $(1, 0.85368, 1.45259)$ & $(1, 1.01573,   0.229231)$    \\
  $E_6/SU(3)\times SU(3)\times SU(2)\times U(1)$ & $(1, 0.771752,  1.33186)$ & $(1, 1.04268,  0.373467)$  \\
   $F_4/SU(3)\times SU(2)\times U(1)$  & $(1, 0.678535,   1.20122)$ & $(1, 1.09057,  0.546045)$ \\
   $G_2/U(2)$   & $(1, 1.67467,   2.05238)$ & $(1, 0.186894,  0.981478)$      \\
  \hline 
       \end{tabular}
        \end{center}}
 \end{theorem}

\markboth{ Stavros Anastassiou and  Ioannis Chrysikos}{The Ricci flow approach to homogeneous Einstein metrics on flag manifolds}
 \section{The method of Poincar\'{e} compactification}
\markboth{ Stavros Anastassiou and  Ioannis Chrysikos}{The Ricci flow approach to homogeneous Einstein metrics on flag manifolds}

Our aim is to calculate the Einstein metrics presented above once again, using the fact that they correspond to fixed points of 
the normalized Ricci flow. As we shall see, the fixed points of this flow are located at infinity and not 
in the finite space. The study of a vector field at infinity is possible by making use of the compactification procedure due to 
Poincar\'{e} (\cite{Poi}).  To make more clear the arguments that follow, we briefly recall the two--dimensional case here, providing also the 
formulas for the three--dimensional case we will need in the subsequent calculations. 
For details on the general case we refer the   reader to \cite{Gon}. 

To study a (polynomial) vector field in a neighborhood of infinity, we introduce a new vector field, defined on a sphere, as follows.  
Let $(x_1,x_2)$ be coordinates on $\mathbb{R}^2$ and $X=P(x_1,x_2)\frac{\partial}{\partial x_1}+Q(x_1,x_2)\frac{\partial}{\partial x_2}$ a polynomial vector field of degree $d$ (that is, 
$d={\rm max}\{{\rm deg}(P),{\rm deg}(Q)\}$). If $(y_1,y_2,y_3)$ 
denote the coordinates on $\mathbb{R}^3$ then we consider $\mathbb{R}^2$ to be the plane of $\mathbb{R}^3$ 
defined as $(y_1,y_2,y_3)=(x_1,x_2,1)$. We consider also the sphere $\mathbb{S}^2=\{y\in \mathbb{R}^3/y_1^2+y_2^2+y_3^2=1\}$, 
which we shall call Poincar\'{e} sphere. This sphere is divided to the northern $(H^+=\{y\in \mathbb{S}^2/y_3>0 \})$ 
and to the southern $(H^-=\{y\in \mathbb{S}^2/y_3<0 \})$ hemispheres, and the equator 
$\mathbb{S}^1=\{y\in \mathbb{S}^2/y_3=0 \}$.
 
 The central projections from $\mathbb{R}^2$ to the Poincar\'{e} sphere are defined as follows:
\begin{center}
$f^+:\mathbb{R}^2\rightarrow \mathbb{S}^2,(x_1,x_2)\mapsto (\frac{x_1}{\Delta (x)},\frac{x_2}{\Delta(x)},\frac{1}{\Delta(x)})$ and \\
$f^-:\mathbb{R}^2\rightarrow \mathbb{S}^2,(x_1,x_2)\mapsto (\frac{-x_1}{\Delta (x)},\frac{-x_2}{\Delta(x)},\frac{-1}{\Delta(x)})$,
\end{center} 
where $\Delta(x)=\sqrt{x_1^2+x_2^2+1}$, and in this way we obtain one vector field on each hemisphere. 
Each one of these vector fields, namely 
\[
\bar {X}(y)=D_xf^+(X(x)), \ y=f^+(x)  \ \ \mbox{and} \ \  \bar{X}(y)=D_xf^-(X(x)), \ y=f^-(x),
\]
 is conjugate to the original vector field. We have thus 
constructed a vector field $\bar{X}$ on $\mathbb{S}^2 \backslash \mathbb{S}^1$ and we want to extend it to $\mathbb{S}^2$. To 
achieve this, we multiply the vector field by the function $\rho (y)=y_3^{d-1}$. We state, without 
proof, the following theorem (\cite{Gon}):
\begin{theorem}
The field $\bar{X}$ can be analytically extended to the whole sphere by multiplication with the factor $y_3^{d-1}$, in such a way 
that the equator is invariant.
\end{theorem} 
\begin{center}
\begin{figure}
\includegraphics[width=0.2\textheight]{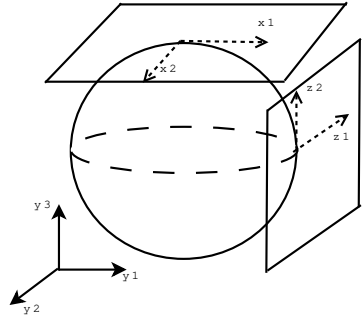}\hspace{2.cm}
\caption{The Poincar\'{e} sphere, $\mathbb{R}^2$ as its tangent space at the north pole and the $U_1$ chart (vertical plane).}
\end{figure}
\end{center}

The vector field defined on the Poincar\'{e} sphere is called the Poincar\'{e} compactification 
of the original vector field $X=(P,Q)$ and is denoted by $p(X)$. Points of the equator correspond 
to the points at infinity of the plane. 

To study the vector field $p(X)$ we make use of six local charts on the Poincar\'{e} sphere, 
given by $U_k=\{y\in \mathbb{S}^2/y_k>0 \}$ and $V_k=\{y\in \mathbb{S}^2/y_k<0 \}$, $k=1,2,3$. 
The local maps for the corresponding charts are given by
\begin{center}
$\phi _k:U_k\rightarrow \mathbb{R}^2$ and $\psi  _k:V_k\rightarrow \mathbb{R}^2$,
\end{center}  
with $\phi_k(y)=-\psi_k(y)=(y_m/y_k,y_n/y_k)$, for $m<n$ and $m,n\neq k$. If we write 
$z=(z_1,z_2)$ for the value of $\phi_k(y)$ or $\psi_k(y)$ then, in any chart, points at infinity 
correspond to $z_2=0$. Note that the meaning of $z$ depends on the chart.

We will now write down the expressions of $p(X)$ in the local charts, for future 
reference. In the chart $(U_1,\phi_1)$ (corresponding in Figure 1 as the $y_1=1$ plane) the expression of the field reads as:
\begin{eqnarray}
\label{u1}
\dot{z}_1&=&z_2^d[-z_1P(\frac{1}{z_2},\frac{z_1}{z_2})+Q(\frac{1}{z_2},\frac{z_1}{z_2})]\\
\dot{z}_2&=&-z_2^{d+1}P(\frac{1}{z_2},\frac{z_1}{z_2}). \nonumber
\end{eqnarray}
The expression in the chart $(U_2,\phi_2)$ is
\begin{eqnarray*}
\dot{z}_1&=&z_2^d[P(\frac{z_1}{z_2},\frac{1}{z_2})-z_1Q(\frac{z_1}{z_2},\frac{1}{z_2})]\\
\dot{z}_2&=&-z_2^{d+1}Q(\frac{z_1}{z_2},\frac{1}{z_2}), \nonumber
\end{eqnarray*}
while for $(U_3,\phi_3)$ is
\begin{eqnarray*}
\dot{z}_1&=&P(z_1,z_2)\\
\dot{z}_2&=&Q(z_1,z_2). \nonumber
\end{eqnarray*}
We omit the expressions of $p(X)$ in the charts $(V_k,\psi_k)$, since 
the coincide with the expressions for $(U_k,\phi_k)$ multiplied by the factor 
$(-1)^{d-1}$ for $k=1,2,3$.

If one is interested to study the global behavior of a vector field 
on $\mathbb{R}^2$, that is, if we are also interested to study the 
vector field in a neighborhood of infinity, then clearly, it is enough 
to work on $H^+\cup \mathbb{S}^1$, which is called the Poincar\'{e} 
disk.

The above described procedure generalizes to every dimension. Since we also need the three--dimensional case 
here we include now the necessary formulas of the compactified vector field, using again $(z_1,z_2,z_3)$ as 
coordinates. If the original vector field is  $X=(P^1,P^2,P^3)$,
 then the equations of the 
compactified field $p(X)$ read as:
\begin{eqnarray*}
\dot{z}_1&=&Q^1(z_1,z_2,z_3)=Q^{1}\\
\dot{z}_2&=&Q^2(z_1,z_2,z_3)=Q^{2}\\
\dot{z}_3&=&Q^3(z_1,z_2,z_3)=Q^{3}
\end{eqnarray*}
where
{\small{ \begin{eqnarray*}
 (Q^1, Q^2, Q^3)&=&   \frac{z_3^d}{(\Delta z)^{d-1}}\Big(-z_1P^1(1/z_3,z_1/z_3,z_2/z_3)+P^2(1/z_3,z_1/z_3,z_2/z_3),  \\
  && -z_2P^1(1/z_3,z_1/z_3,z_2/z_3)+P^3(1/z_3,z_1/z_3,z_2/z_3),-z_3P^1(1/z_3,z_1/z_3,z_2/z_3)\Big),  \\
 (Q^1, Q^2, Q^3)&=&\frac{z_3^d}{(\Delta z)^{d-1}}\Big(-z_1P^2(z_1/z_3,1/z_3,z_2/z_3)+P^1(z_1/z_3,1/z_3,z_2/z_3), \\
 && -z_2P^2(z_1/z_3,1/z_3,z_2/z_3)+P^3(z_1/z_3,1/z_3,z_2/z_3),-z_3P^2(z_1/z_3,1/z_3,z_2/z_3)\Big), \\
 (Q^1, Q^2, Q^3)&=& \frac{z_3^d}{(\Delta z)^{d-1}}\Big(-z_1P^3(z_1/z_3,z_2/z_3,1/z_3)+P^1(z_1/z_3,z_2/z_3,1/z_3),\\ && -z_2P^3(z_1/z_3,z_2/z_3,1/z_3)+P^2(z_1/z_3,z_2/z_3,1/z_3),-z_3P^3(z_1/z_3,z_2/z_3,1/z_3)\Big),\\
\end{eqnarray*}    
in   $U_1$, $U_2$ and $U_3$, respectively.
The expression for $p(X)$ in $U_4$ is  
\[
z_3^{d+1}(P^1(z_1,z_2,z_3),P^2(z_1,z_2,z_3),P^3(z_1,z_2,z_3)),
\]
 while 
the expressions at the local charts $V_i$ can be obtained by those at the 
charts $U_i$ multiplied by the term $(-1)^{d-1}$. By rescaling the time 
variable, we usually omit the factor $1/(\Delta z)^{d-1}$.

This procedure will be now used, with the aim of analyzing the global behavior of the Ricci flow 
equations presented in the previous section.

 \markboth{ Stavros Anastassiou and  Ioannis Chrysikos}{The Ricci flow approach to homogeneous Einstein metrics on flag manifolds}
 \section{The global behaviour of the normalized Ricci flow}
 \markboth{ Stavros Anastassiou and  Ioannis Chrysikos}{The Ricci flow approach to homogeneous Einstein metrics on flag manifolds}
 
\subsection{Dynamics of the   Ricci flow on flag manifolds with $\fr{m}=\fr{m}_1\oplus\fr{m}_2$}
In this section we  study the global behavior of the normalized Ricci flow 
equation for a flag manifold $M=G/K$ with two isotropy summands,  namely  system (\ref{nrf2}).  This system reduces to:
\begin{eqnarray}
\dot{x}_1&=& \displaystyle\frac{8d_2^2x_1^2+2(2d_1+d_2)(d1+4d_2)x_1x_2-d_2(3d_1+2d_2)x_2^2}{2(d_1+d_2)(d_1+4d_2)x_1x_2}\nonumber\\
\dot{x}_2&=& \displaystyle\frac{(4d_2x_1+d_1x_2)(4d_2x_1+d_1(2x_1+x_2))}{2(d_1+d_2)(d_1+4d_2)x_1^2}.\nonumber
\end{eqnarray}
To apply the Poincar\'{e} compactification method we multiply this system with the factor 
$2(d_1+d_2)(d_1+4d_2)x_1^2x_2$.  We remark here that this multiplication will only change the time 
parametrization of the orbits and not the structure of the phase portrait we wish to determine. 
We arrive thus to the following two--dimensional system:
\begin{equation}\label{ricci2d}
\left. \begin{tabular}{ll}
$\dot{x}_1=$ & $8d_2^2x_1^2+2(2d_1+d_2)(d1+4d_2)x_1x_2-d_2(3d_1+2d_2)x_2^2$, \\
$\dot{x}_2=$ & $(4d_2x_1+d_1x_2)(4d_2x_1+d_1(2x_1+x_2))$,  
\end{tabular}\right\}
\end{equation}
the behavior of which, at the first quadrant of the plane, is to be determined. We consider 
$d_1,d_2>0$ as free parameters of the system.

We begin our study with the following lemma.
\begin{lemma}
\label{finite2d}
System (\ref{ricci2d}) possesses a single fixed point at the origin. Except from the coordinate 
axes, the straight lines 
$\gamma _1(t)=(\frac{1}{2}t,t),\gamma _2(t)=(\frac{d_1+2d_2}{4d_2}t,t)$, 
remain invariant under its flow.
\end{lemma} 
We omit the proof of Lemma \ref{finite2d}, since it is quite straightforward. 

From the dynamical viewpoint, the invariant axes and the lines 
$\gamma_1,\gamma_2$ consist the separatrices of the parabolic sectors 
of the fixed point located at the origin (at least the separatrices 
in the first quadrant). From the Ricci flow viewpoint, the invariance 
of the axes reflect the fact that no semi--Riemannian or degenerate 
metric can evolve to a Riemannian metric via the Ricci flow, while, as 
we shall see, the other two straight lines are related to the Einstein 
metrics that we are about to compute.

We proceed now to the study of system (\ref{ricci2d}) at infinity. Since we are only interested 
in positive values of $x_1$ and $x_2$ we restrict ourselves to the $(U_1,\phi_1)$ chart. 
Applying (\ref{u1}) and using the notation of the previous section we easily confirm that equations (\ref{ricci2d}), 
in this chart, take the following form:
\begin{eqnarray}
\dot{z}_1&=&(d_1+d_2)(z_1-2)z_1(2d_2(z_1-2)+z_1d_1),\nonumber\\
\dot{z}_2&=&(-4d_1^2z_1+3d_1d_2(z_1-6)z_1+2d_2^2(-4+(z_1-4)z_1))z_2.\nonumber
\end{eqnarray}
By setting $z_2=0$ we calculate that there are only two fixed points 
(except from the origin), namely $(2,0)$ and $(\frac{4d_2}{d_1+2d_2},0)$. 
The first one is a repelling node and the other 
one an attracting node. The two invariant lines previously calculated converge to these 
fixed points, and we are thus able to draw the global phase portrait of the system 
in Figure 2.
\begin{center}
\begin{figure}
\includegraphics[width=0.3\textheight]{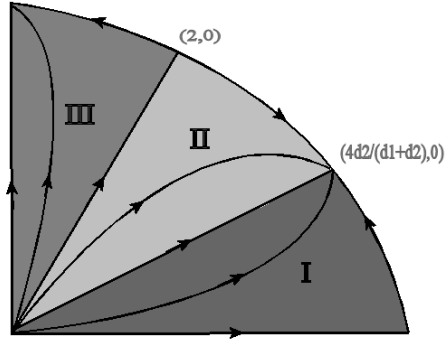} \hspace{1.cm}
\caption{System (\ref{ricci2d}) in the first quadrant of the plane.}
\end{figure}
\end{center}
\begin{prop}
The global phase portrait of system (\ref{ricci2d}) is topologically equivalent to the one 
depicted in Figure (2), when we restrict our attention to the  first quadrantant 
of the plane. 
\end{prop}

We can now draw conclusions about the existence of Einstein metrics. Recalling that the chart $(U_1,\phi_1)$ 
corresponds to the $y_1=1$ plane (see the section on the Poicar\'{e} compactification), we consider the metrics whose coefficients are determined by the fixed points calculated above:
\[
g_1= 1\cdot ( \ , \ )|_{\fr{m}_1}+2\cdot ( \ , \ )|_{\fr{m}_2}, \quad \mbox{and} \quad g_2 =1\cdot ( \ , \ )|_{\fr{m}_1}+\frac{4d_2}{d_1+2d_2}\cdot ( \ , \ )|_{\fr{m}_2}.
\]
We obtain the following:
\begin{theorem}
Let $M=G/K$ be a generalized flag manifold of a compact simple Lie group $G$, with $\fr{m}=\fr{m}_1\oplus\fr{m}_{2}$. The normalized Ricci flow, on the space of invariant Riemannian metrics, possesses exactly two fixed 
points at infinity, one of them being an attracting node while the other one is a repelling one. These fixed points correspond to 
the  two $G$--invariant Einstein metrics of $M$, namely $g_1,g_2$ defined above. The 
first metric is the (unique) K\"{a}hler--Einstein metric on $M$ while the second one 
is non K\"{a}hler. Every initial metric belonging to region I and II or the line $\gamma_2$ (see Figure (2)) tends to 
the $g_2$ metric, while the metrics belonging to the line $\gamma_1$ tend to the K\"{a}hler--Einstein metric, 
under the normalized Ricci flow. 
\end{theorem}

\subsection{Dynamics of the normalized Ricci flow on flag manifolds $M=G/K$ of Type I}
In this section we  study the global behavior of the normalized Ricci flow 
equation for a flag manifold $M=G/K$ with three isotropy summands and Type I,  namely the  system (\ref{nrf3}).  This system is not polynomial, but after multiplication with the (positive, in the first octand) factor $2d_1d_2d_3(d_1+d_2+d_3)(d_1+4d_2+9d_3)x_1^2x_2x_3$ it reduces to:
{\footnotesize{
\begin{eqnarray}\label{ricci3d}
\dot{x}_{1}  &=&  d_2 x_1 \Big(2 d_3 x_1 \big((d_1 + d_2) (d_2 + d_3) x_1^2 + 
      d_1 (d_1 + 4 d_2 + 9 d_3) x_1 x_2 - (d_1 + d_2) (2 d_1 + d_2 + 
         d_3) x_2^2\big) \nonumber \\
         && + \big(-4 d_1 (-2 d_2^2 + d_1 d_3 - 5 d_2 d_3) x_1^2 + 
      2 d_1 (2 d_1 + d_2 + d_3) (d_1 + 4 d_2 + 9 d_3) x_1 x_2 \nonumber \\
      && - (3 d_1 +  2 (d_2 + d_3)) (-d_2 d_3 + d_1 (d_2 + 2 d_3)) x_2^2\big) x_3  
          -    2 (d_1 + d_2) d_3 (2 d_1 + d_2 + d_3) x_1 x_3^2\Big), \nonumber \\
\dot{x}_{2}&=& d_1 x_2 \Big(2 d_3 x_1 \big(-(d_1 + d_2) (d_1 + 2 d_2 + d_3) x_1^2 + 
      d_2 (d_1 + 4 d_2 + 9 d_3) x_1 x_2 + (d_1 + d_2) (d_1 + 
         d_3) x_2^2\big) \nonumber \\
         && + \big(-4 (d_1 + 2 d_2 + d_3) (-2 d_2^2 + d_1 d_3 - 
         5 d_2 d_3) x_1^2 + 
      2 d_1 d_2 (d_1 + 4 d_2 + 9 d_3) x_1 x_2 \nonumber \\
      && + (d_1 + d_3) (-d_2 d_3 + 
         d_1 (d_2 + 2 d_3)) x_2^2\big) x_3 - 
   2 (d_1 + d_2) d_3 (d_1 + 2 d_2 + d_3) x_1 x_3^2), \nonumber \\
   \dot{x}_{3}&=& d_1 d_2 x_3 \Big(-2 (d_1 + d_2 + 
      2 d_3) x_1 \big((d_1 + d_2) x_1^2 - (d_1 + 4 d_2 + 9 d_3) x_1 x_2 + (d_1 + 
         d_2) x_2^2\big) \nonumber \\
         &&  + \big(4 (2 d_2^2 - d_1 d_3 + 5 d_2 d_3) x_1^2 + 
      2 d_1 (d_1 + 4 d_2 + 9 d_3) x_1 x_2 \nonumber \\
      && + (d_2 d_3 - 
         d_1 (d_2 + 2 d_3)) x_2^2\big) x_3 + 2 (d_1 + d_2)^2 x_1 x_3^2\Big).  
\end{eqnarray}}}
The analysis of a 3--dimensional dynamical system, like the one above, is of course a quite challenging task, 
but since we are  interested here only in the most elementary of its properties we may proceed as follows.

We once again remark that we study this system in the region $P=\{(x_1,x_2,x_3)\in\mathbb{R}^3/x_1,x_2,x_3>0\}$, since 
we are interested in Riemannian metrics. In this region, and in analogy with Lemma (\ref{finite2d}),  we have the following:
\begin{lemma}
\label{finite3d}
The coordinate planes, along with the straight line $\rho(t)=(t,2t,3t)$, remain invariant under the 
flow the system (\ref{ricci3d}) defines. Moreover, this system possesses no fixed points in the region $P$, 
for any value of the parameters $(d_1,d_2,d_3)$ reported in Table 4.
\end{lemma}
Translated into the language of Ricci flow, the above lemma ensures that the normalized Ricci flow possesses 
no singularities in finite region. Moreover, the invariance of the coordinate planes prohibits the evolution of 
a degenerate, or semi--Riemannian metric to a Riemannian one, while the invariant straight line is once again 
related with the existence of an Einstein metric we are about to obtain.

Since there no singularities in the finite region, we now turn our attention to the study if system (\ref{ricci3d}) 
at infinity. We apply the Poincar\'{e} compactification procedure previously described and arrive to the following system, written in 
the $(U_1,\phi_1)$ chart:
{\footnotesize{
 \begin{eqnarray}
\dot{z}_{1}&=&  (d_1  + d_2 + d_3) z_1 \Big(2 d_2^2 d_3 \big(-1 - z_1^2 (-1 + z_2) + z_2^2\big) + 
   d_1^2 \big(d_2 (-2 + z_1) z_1 z_2 + \nonumber \\
   && 2 d_3 (-1 + z_1^2 - z_2) (1 + z_2)\big)  
    +  d_1 d_2 \big(2 d_2 (-2 + z_1)^2 z_2 + 
      d_3 (-4 + 20 z_2 + z_1 (4 z_1 + 3 (-6 + z_1) z_2))\big)\Big), \nonumber \\
\dot{z}_2 &=&  -2 d_2 (d_1 + d_2 + d_3) z_2 \Big(d_2 d_3 (-1 + z_1^2 - z_2) (-1 + z_2) + 
   d_1^2 (1 + z_1^2 + z_1 (-1 + z_2) - z_2^2) \nonumber \\
   &&  +  d_1 \big(d_3 - d_2 (-1 + z_2) (1 + (-4 + z_1) z_1 + z_2) - 
      d_3 (z_1 (9 + z_1) + z_1 (-9 + 2 z_1) z_2 + z_2^2)\big)\Big), \nonumber\\
\dot{z}_{3}&=& -d_2 \Big(2 d_3 \big((d_1 + d_2) (d_2 + d_3) + 
      d_1 (d_1 + 4 d_2 + 9 d_3) z_1 - (d_1 + d_2) (2 d_1 + d_2 + 
         d_3) z_1^2\big) \nonumber \\
         && + \big(4 d_1 (2 d_2^2 - d_1 d_3 + 5 d_2 d_3) + 
      2 d_1 (2 d_1 + d_2 + d_3) (d_1 + 4 d_2 + 9 d_3) z_1 - \nonumber \\
      && (3 d_1 + 
         2 (d_2 + d_3)) (-d_2 d_3 + d_1 (d_2 + 2 d_3)) z_1^2\big) z_2 - 
   2 (d_1 + d_2) d_3 (2 d_1 + d_2 + d_3) z_2^2\Big) z_3. \nonumber
\end{eqnarray}}}
The behavior at infinity is governed by the system obtained after setting $z_3=0$ to the previous equations. It reads as:
{\small{
\begin{eqnarray}\label{22}
\dot{z}_1&=& (d_1 + d_2 + 
   d_3) z_1 \Big(2 (d_1 + d_2)^2 d_3 (-1 + z_1^2) + \big(8 d_1 d_2^2 - 
      4 d_1 (d_1 - 5 d_2) d_3 - 
      2 d_1 d_2 (d_1 + 4 d_2 + 9 d_3) z_1 \nonumber \\
      && + (d_1 + 2 d_2) (-d_2 d_3 + 
         d_1 (d_2 + 2 d_3)) z_1^2\big) z_2 + 2 (-d_1^2 + d_2^2) d_3 z_2^2\Big),  \nonumber \\
\dot{z}_2&=&  -2 d_2 (d_1 + d_2 + d_3) z_2 \Big(d_2 d_3 (-1 + z_1^2 - z_2) (-1 + z_2) + 
   d_1^2 (1 + z_1^2 + z_1 (-1 + z_2) - z_2^2) \nonumber \\
   &&  + 
   d_1 \big(d_3 - d_2 (-1 + z_2) (1 + (-4 + z_1) z_1 + z_2) - 
      d_3 (z_1 (9 + z_1) + z_1 (-9 + 2 z_1) z_2 + z_2^2)\big)\Big).
 \end{eqnarray}
}}
It is not difficult to verify that the system above possesses always a singularity, located at $(2,3)$, which 
is a repelling node. Since it is complex enough to prevail the analytical study of other fixed points, we substitute 
the values of the dimensions $d_1,d_2,d_3$ reported in Table 4, and study each case seperately, given the corresponding form of system (\ref{22}).  In any case we compute two more fixed points which are saddles.
\begin{itemize}
\item{$G_2/U(2)$ 
 \begin{eqnarray*}
\dot{z}_1&=&10 z_1 (288 (-1 + z_1^2) + (512 - 768 z_1 + 256 z_1^2) z_2 - 96 z_2^2)\\
\dot{z}_2&=&640 z_2 (-3 + 12 z_1 + 2 (-6 + z_1) z_1 z_2 + 3 z_2^2)
\end{eqnarray*}
 Equillibria: $(2,3)$, $(0.186894,0.981478)$ and $(1.67467,2.05238)$.}
\item{$E_6/SU(3)\times SU(3)\times SU(2)\times U(1)$   
\begin{eqnarray*}
\dot{z}_1&=&58 z_1 (23328 (-1 + z_1^2) + (124416 - 186624 z_1 + 62208 z_1^2) z_2 - 
   7776 z_2^2)\\
\dot{z}_2&=&902016 z_2 (-5 - 4 (-3 + z_1) z_1 + 2 (-6 + z_1) z_1 z_2 + 5 z_2^2)
\end{eqnarray*}
Equillibria: $(2, 3)$,  $(0.771752,1.33186)$ and 
$(1.04268,0.373467)$.}
\item{$E_7/SU(5)\times SU(3)\times U(1)$ 
\begin{eqnarray*}
\dot{z}_1&=&4233600 z_1 (3 (-1 + z_1^2) + 7 (2 + (-3 + z_1) z_1) z_2 - z_2^2)\\
\dot{z}_2&=&2116800 z_2 (-17 + 42 z_1 - 13 z_1^2 + 7 (-6 + z_1) z_1 z_2 + 17 z_2^2)
\end{eqnarray*}
 Equillibria: $(2,3)$, $(0.733552,1.27681)$ and $(1.06029,0.443559)$.}
\item{$E_7/SU(6)\times SU(2)\times U(1)$ 
\begin{eqnarray*}
\dot{z}_1&=&2030400 z_1 (3 (-1 + z_1^2) + 12 (2 + (-3 + z_1) z_1) z_2 - z_2^2)\\
\dot{z}_2&=&4060800 z_2 (-8 + 18 z_1 - 7 z_1^2 + 3 (-6 + z_1) z_1 z_2 + 8 z_2^2)
\end{eqnarray*}
 Equillibria: $(2,3)$, $(0.85368,1.45259)$ and $(1.01573,0.229231)$.}
\item{$E_8/E_6\times SU(2)\times U(1)$ 
\begin{eqnarray*}
\dot{z}_1&=&11617344 z_1 (3 (-1 + z_1^2) + 20 (2 + (-3 + z_1) z_1) z_2 - z_2^2)\\
\dot{z}_2&=&23234688 z_2 (-14 + 30 z_1 - 13 z_1^2 + 5 (-6 + z_1) z_1 z_2 + 14 z_2^2)
\end{eqnarray*} 
 Equillibria: $(2,3)$, $(0.914286,1.54198)$ and $(1.0049,0.129681)$.}
\item{$E_8/SU(8)\times U(1)$ 
\begin{eqnarray*}
\dot{z}_1&=&18464768 z_1 (9 (-1 + z_1^2) + 20 (2 + (-3 + z_1) z_1) z_2 - 3 z_2^2)\\
\dot{z}_2&=&36929536 z_2 (-30 z_1 (-1 + z_2) + z_1^2 (-9 + 5 z_2) + 12 (-1 + z_2^2))
\end{eqnarray*} 
 Equillibria: $(2,3)$, $(0.177586,1.25432)$ and  $(1.06853,0.473177)$.}
\item{$F_4/SU(3)\times SU(2)\times U(1)$ 
\begin{eqnarray*}
\dot{z}_1&=&138240 z_1 (3 (-1 + z_1^2) + 6 (2 + (-3 + z_1) z_1) z_2 - z_2^2)\\
\dot{z}_2&=&138240 z_2 (-7 + 18 z_1 - 5 z_1^2 + 3 (-6 + z_1) z_1 z_2 + 7 z_2^2)
\end{eqnarray*} 
 Equillibria: $(2,3)$, $(0.678535,1.20122)$ and $(1.09057,0546045)$.}
\end{itemize}

Since the $(U_1,\phi_1)$ chart corresponds to the $y_1=1$ plane, we consider the metrics whose coefficients 
are given by $(1,a,b)$, where $a,b$ are equal to the coordinates of the fixed points just calculated. These metrics are 
invariant Einstein metrics and the one with coefficients $(1,2,3)$ is the unique K\"{a}hler--Einstein which admits $M$. We have 
thus proved the following:
\begin{theorem}
Let $M=G/K$ be a generalized flag manifold of a compact simple Lie group 
$G$ of Type I. The normalized Ricci flow, on the space of invariant Riemmanian metrics on $M$, possesses exactly 
three singularities at infinity. The point $(2,3)$ is a reppeling node while the other two are saddle points. These 
fixed points correspond to the three (up to scale) $G$-invariant Einstein matrics which admits  $M$.
\end{theorem}

\subsection{Conclusions}
 We studied in this paper the behavior of the normalized Ricci flow on flag manifolds with two or three 
isotropy summands. The Ricci flow equation reduces to a system of two, correspondingly   three,
ordinary differential equations.

In the case of two isotropy summands we were able to completely determine the system's global 
phase portrait (at the first quadratant of the plane, since we are interested in Riemmanian metrics), 
using standard techniques of Dynamical Systems theory.  In the case of three isotropy 
summands the problem becomes (as usual) more complicated, making the presentation of the 
corresponding complete phase portrait impossible. In any case, we were able to calculate 
explicitly the invariant Einstein metrics that exist, as the fixed points of the normalized Ricci flow at infinity. Morever, we can point out the 
K\"{a}hler--Einstein metric as the fixed point having no stable eigendirections.

Obtaining the invariant Einstein metrics using the Ricci flow does not only provides their explicit number 
and form, but also eluminates the asymptotic behavior of $G$-invariant Riemmanian metrics.   This unified method may prove to be useful in obtaining 
more general results about the existence and the classification problem   of homogeneous Einstein metrics on compact homogeneous spaces.



\begin{thebibliography}{50}
 \bibitem {AA} 
 D. V. Alekseevsky and A. Arvanitoyeorgos: 
{\it Riemannian flag manifolds with homogeneous geodesics},
 Trans. Amer. Math. Soc. 359 (8)  (2007)  3769--3789. 
 
 
   \bibitem {AP} 
 D. V. Alekseevsky and A. M. Perelomov: 
{\it Invariant K\"ahler-Einstein metrics on compact homogeneous spaces}, 
Funct. Anal. Appl. 20 (3)  (1986)  171--182. 

\bibitem {Arv}  
A. Arvanitoyeorgos: 
{\it An Introduction to Lie Groups and the Geometry of Homogeneous Spaces}, 
Amer. Math. Soc, Student Math. Library, Vol. 22, 2003. 

 \bibitem  {Chry1} 
 A. Arvanitoyeorgos and I. Chrysikos:
 {\it Motion of charged particles and homogeneous geodesics 
 in K\"ahler $C$-spaces with two isotropy summands}, 
 Tokyo J. Math. (32) 2 (2009)  1-14.
 
 
  \bibitem {Chry2} 
A. Arvanitoyeorgos and I. Chrysikos:
{\it Invariant Einstein metrics on generalized 
flag manifolds with two isotropy summands}, 
to appear in J. Aust. Math. Soc. (2010). 
 
 
   \bibitem {Chry3} 
A. Arvanitoyeorgos and I. Chrysikos:
{\it Invariant Einstein metrics on generalized flag manifolds with four isotropy summands}, 
 Ann. Glob. Anal. Geom. 37 (2010) 185-219.
 
  \bibitem {Be}
 A. L. Besse:
  {\it Einstein Manifolds},
  Springer-Verlag, Berlin, 1986.
  
    \bibitem {Bohm}
   C. B\"ohm and B. Wilking:
   {\it Nonnegatively curved manifolds with finite 
   fundamental groups admits metrics with positive Ricci curvature},
    Geom. Funct. Anal. 17 (2007)   665-681.
    
     \bibitem{Borel} 
      A. Borel and F. Hirzebruch: 
      {\it Characteristic  classes and homogeneous spaces I}, 
       Amer. J. Math. 80  (1958)  458--538.
       
        \bibitem {CaCh}
  H-D. Chao and B. Chow:
  {\it Recent developments on the Ricci flow},
  Bull. Amer. Soc. 36 (1999) 59-74.  
  
   \bibitem {Do} 
    R. Dohira: 
    {\it Geodesics in reductive homogeneous spaces}, 
     Tsukuba J. Math. 19 (1) (1995) 233-243. 
     
     \bibitem{Kim}  
 M. Kimura: 
 {\it Homogeneous Einstein metrics on certain K\"ahler C-spaces},
  Adv. Stud. Pure Math.  18-I  (1990)  303--320. 
  
     \bibitem {Ker}
  W. Dickinson and M. Kerr:
  {\it The geometry of compact homogeneous spaces with two isotropy summands}, 
  Ann. Glob. Anal. Geom. 34 (2008)   329–-350.
  
  \bibitem {Glick}
  D. Glickenstein and T. L. Payne:
  {\it Ricci flow on three-dimensional, unimodular metric Lie algebras}, 
   preprint  (2009),  arXiv:0909.0938v1.
   
  \bibitem {Gon}
E.A.V. Gonzales: {\it Generic properties of polynomial vector fields at infinity}, 
Trans. Amer. Math. Soc. 143:201-222, 1969.  

  \bibitem {Grama}
L. Grama and R. M. Martins:
{\it The Ricci flow of left-invariant metrics on full 
flag manifold $SU(3)/T$ from a dynamical systems point of view},
Bull. Sci. Math. 133 (2009) 463-469.
, 
\bibitem {Hel}
 S. Helgason:
  {\it Differential Geometry, Lie Groups and Symmetric Spaces},
   Academic Press, New York 1978
   
\bibitem {Ham1}
R. S. Hamilton:
{\it Three-manifolds with positive Ricci curvature}, 
J. Differential Geom. 17 (1982) 255-306.

    \bibitem {Ham2}
R. S. Hamilton:
{\it Four-manifolds with positive  curvature operator}, 
J. Differential Geom. 24 (1986) 153-179.

      \bibitem {SP}
  J-S. Park and Y. Sakane:
  {\it Invariant Einstein metrics on certain homogeneous spaces},  
  Tokyo J. Math.  20   (1) (1997) 51--61.
  
  \bibitem{Poi} 
  H. Poincar\'{e}: {\it Sur l' integration des \'{e}quations diff\'{e}rentielles du premier ordre et du premier degr\'{e} I}, Rendiconti del circolo matematico di Palermo,  5,  (1891) 161-191.
  \bibitem{Sak}
  Y. Sakane
  {\it Homogeneous Einstein metrics on principal circle bundles II},
  in Differential Geometry, Proceedings of the symposioum in honour 
  of Professor Su Buchin on his 90th birthday (Editors: C. H. Gu, H. S. Hu, Y. L. Xin), World Scientific Publishing, (1993)   177--186.
  
       \bibitem {Wa2}
  M. Wang and W. Ziller: 
  {\it Existence and non-excistence of homogeneous Einstein metrics}, 
  Invent.~Math.~84 (1986)  177--194.
  \end{thebibliography}
\end{document}